# WEIGHTED POINCARÉ-TYPE INEQUALITIES FOR CAUCHY AND OTHER CONVEX MEASURES[1]

BY SERGEY G. BOBKOV AND MICHEL LEDOUX

*University of Minnesota and Université Paul-Sabatier*

Brascamp–Lieb-type, weighted Poincaré-type and related analytic inequalities are studied for multidimensional Cauchy distributions and more general $\kappa$-concave probability measures (in the hierarchy of convex measures). In analogy with the limiting (infinite-dimensional log-concave) Gaussian model, the weighted inequalities fully describe the measure concentration and large deviation properties of this family of measures. Cheeger-type isoperimetric inequalities are investigated similarly, giving rise to a common weight in the class of concave probability measures under consideration.

**1. Introduction.** The aim of these notes is to study some aspects of the high-dimensional analysis, such as measure concentration and isoperimetric properties, of families of convex probability measures through inequalities of the type of Brascamp–Lieb, Poincaré, Cheeger and logarithmic Sobolev, with weight. Although such inequalities may be considered in different contexts and settings, we restrict ourselves to the Euclidean space $\mathbb{R}^n$. The multidimensional Cauchy distribution is the prototype model in the family of convex measures under investigation, and may be analyzed in analogy with the Gaussian model for the class of log-concave measures. Under a proper scaling, the Gaussian model actually appears as the limiting case of the (finite-dimensional) Cauchy model.

A Borel probability measure $\mu$ on $\mathbb{R}^n$ is said to satisfy a weighted Poincaré-type inequality with weight function $w^2$ (where $w$ is a fixed nonnegative, Borel measurable function), if for any bounded smooth function $g$ on $\mathbb{R}^n$

Received April 2007; revised April 2008.
[1]Supported in part by the NSF Grant DMS-07-06866.

with gradient $\nabla g$,

$$\text{(1.1)} \qquad \text{Var}_\mu(g) \leq \int |\nabla g|^2 w^2 \, d\mu.$$

As usual, $\text{Var}_\mu(g) = \int g^2 \, d\mu - (\int g \, d\mu)^2$ stands for the variance of $g$ under $\mu$. Throughout this paper, we state the various functional inequalities for smooth and bounded, or square integrable, functions. The inequalities may then be classically extended to the class of all locally Lipschitz functions (Poincaré or similar inequalities are then understood in the following sense: if the right-hand side is finite, then the function is square-integrable, and the inequality holds true).

As a classical example, the standard Gaussian measure $\mu = \gamma_n$ with density $d\gamma_n(x)/dx = (2\pi)^{-n/2} e^{-|x|^2/2}$ with respect to Lebesgue measure on $\mathbb{R}^n$ satisfies (1.1) with $w = 1$ (cf., e.g., [17]). The dimension-free character of this inequality allows one to consider important properties of infinite-dimensional Gaussian measures in terms of general Gaussian processes. An attempt to extend this example to other probability distributions may lead to other forms of the usual Poincaré-type inequalities, including those that admit nonconstant weights. Of particular interest is the family of the generalized Cauchy distributions $\nu_\beta$ on $\mathbb{R}^n$, which have densities

$$\frac{d\nu_\beta(x)}{dx} = \frac{1}{Z}(1 + |x|^2)^{-\beta}$$

with parameter $\beta > n/2$ and a normalizing constant $Z$, depending on $n$ and $\beta$. These densities play an important role in different mathematical problems. For example, they appear to be the extremal functions for the classical Sobolev inequalities in $\mathbb{R}^n$. Under a different name (as a finite-mass Barenblatt profile), they are also used to represent a stationary attractor of the nonlinear diffusion equation $\partial u/\partial t = \Delta(u^m) + \text{div}(xu)$.

Although the family $\nu_\beta$ is of interest in itself, it may be considered in particular as a natural "pre-Gaussian" model, where the Gaussian case appears in the limit as $\beta \to \infty$ (after proper rescaling of the coordinates). This view inspires us to look for general geometric and analytic properties of the Cauchy distributions, which would describe some finite-dimensional analogues of the infinite-dimensional Gaussian model as well as recover a number of known results for the Gaussian measure $\gamma_n$, such as (1.1), in the limit as $\beta \to \infty$. To this aim we prove in particular that each $\nu_\beta$ with $\beta \geq n$ satisfies the weighted Poincaré-type inequality

$$\text{(1.2)} \qquad \text{Var}_{\nu_\beta}(g) \leq \frac{C}{\beta} \int |\nabla g(x)|^2 (1 + |x|^2) \, d\nu_\beta(x)$$

for all bounded smooth functions $g$ on $\mathbb{R}^n$ with a numerical constant $C$. Simple test functions suggest that, up to a constant, the weight function is chosen correctly. Inequality (1.2) can be used to study large deviations of Lipschitz functionals $g$ under $\nu_\beta$. While all exponential moments might be infinite (like for the Euclidean norm), we will see that the tails $\nu_\beta(|g - \int g \, d\nu_\beta| \geq t)$ admit an exponential bound on a long interval with length proportional to $\beta$. Actually, on smaller intervals the tails are sub-Gaussian, and this leads to a certain generalization of the concentration phenomenon for the Gaussian measure. From (1.2) one may also derive a reverse weighted form

$$(1.3) \qquad \inf_{c \in \mathbb{R}} \int \frac{|g-c|^2}{1+|x|^2} \, d\nu_\beta(x) \leq C_\beta \int |\nabla g|^2 \, d\nu_\beta.$$

Using symmetrization and a careful treatment of an associated one-dimensional problem, a similar imbedding inequality was recently obtained by Denzler and McCann [12] and applied there to study the rate of convergence to the stationary attractor of the diffusion equation.

Returning to the Gaussian Poincaré-type inequality, in the mid 1970s, Brascamp and Lieb [10] proposed a remarkable extension to the class of probability measures $\mu$ on $\mathbb{R}^n$ with log-concave densities

$$(1.4) \qquad p(x) = e^{-W(x)}.$$

Namely, let $W$ be a twice continuously differentiable convex function with positive second derivative $W''$ in the matrix sense (so that the inverse matrix $W''^{-1}$ exists and is continuous). Then, for all bounded smooth functions $g$ on $\mathbb{R}^n$,

$$(1.5) \qquad \mathrm{Var}_\mu(g) \leq \int \langle W''^{-1} \nabla g, \nabla g \rangle \, d\mu.$$

The formula (1.4) describes the class of the so-called log-concave probability measures, which are known to satisfy also the usual Poincaré-type inequality, that is, (1.1) with $w = C(\mu)$, compare [4]. For this class, (1.1) and (1.5) are in general not comparable. When however $W'' \geq c \, \mathrm{Id}$ for some $c > 0$ in the sense of symmetric matrices, (1.5) yields

$$(1.6) \qquad \mathrm{Var}_\mu(g) \leq \frac{1}{c} \int |\nabla g|^2 \, d\mu.$$

In full analogy with the Gaussian case, we examine here the natural extension of the family of Cauchy distributions to the so-called class of $\kappa$-concave measures. A Radon probability measure $\lambda$ on a locally convex space $L$ is called $\kappa$-concave, where $-\infty \leq \kappa \leq +\infty$, if it satisfies

$$(1.7) \qquad \lambda_*(tA + (1-t)B) \geq [t\lambda(A)^\kappa + (1-t)\lambda(B)^\kappa]^{1/\kappa}$$

for all $t \in (0,1)$ and for all Borel measurable sets $A, B \subset L$ with positive measure, where $\lambda_*$ stands for the inner measure associated to $\lambda$. When $\kappa = 0$, the right-hand side of (1.7) is understood as $\lambda(A)^t \lambda(B)^{1-t}$ and then we arrive at the notion of a log-concave measure, previously considered by Prékopa [22, 23] and Leindler [18]. When $\kappa = -\infty$, the right-hand side is understood as $\min\{\lambda(A), \lambda(B)\}$. The inequality (1.7) is getting stronger as the parameter $\kappa$ is increasing, so the case $\kappa = -\infty$ describes the largest class whose members are called convex or hyperbolic probability measures. The family of probability measures satisfying the Brunn–Minkowski-type inequality (1.7) was introduced and studied by Borell [8, 9]. In particular, he gave the following characterization of such measures: If $L$ has finite dimension, any $\kappa$-concave probability measure is supported on some convex set $\Omega \subset L$ and is absolutely continuous with respect to Lebesgue measure on $\Omega$. Necessarily, $\kappa \leq 1/\dim(\Omega)$, so if $\lambda$ is not a delta measure, we have $\kappa \leq 1$. If $L$ has dimension $n$ and $\lambda$ is absolutely continuous with respect to Lebesgue measure on $L$, the necessary and sufficient condition that $\lambda$ satisfies (1.7) is that it is supported on some open convex set $\Omega \subset L$, where it has a positive density $p$ such that, for all $t \in (0, 1)$,

$$(1.8) \quad p(tx + (1-t)y) \geq [tp(x)^{\kappa_n} + (1-t)p(y)^{\kappa_n}]^{1/\kappa_n}, \qquad x, y \in \Omega,$$

where $\kappa_n = \kappa/(1 - n\kappa)$ (necessarily $\kappa \leq 1/n$). If $L$ has infinite dimension, then $\lambda$ is $\kappa$-concave if and only if all finite-dimensional projections of $\lambda$ are $\kappa$-concave.

This description may be applied to the generalized Cauchy distribution $\nu_\beta$, in which case $\kappa < 0$, so that the inequality (1.8) turns into

$$(1 + |tx + (1-t)y|^2)^{-\beta \kappa_n} \leq t(1 + |x|^2)^{-\beta \kappa_n} + (1-t)(1 + |y|^2)^{-\beta \kappa_n}.$$

Note that $r = -\beta \kappa_n > 0$ and that the function $z \to (1 + |z|^2)^r$ with $r > 0$ is convex on $\mathbb{R}^n$ if and only if $r \geq 1/2$. Hence, an optimal value of $\kappa$ corresponds to $-\beta \kappa_n = 1/2$. As a conclusion, if $\beta = (n+d)/2$, $d > 0$, the generalized Cauchy measure $\nu_\beta$ is $\kappa$-concave with the optimal value

$$\kappa = -\frac{1}{d} = -\frac{1}{2\beta - n}.$$

For example, the standard Cauchy measure on the line is $\kappa$-concave with $\kappa = -1$. Thus, the characteristics $d$ or $\kappa$, which are directly responsible for convexity properties, may be taken equivalently as the main parameter for the generalized Cauchy measures.

More generally, one may consider probability measures $\mu$ with densities with respect to Lebesgue measure

$$(1.9) \qquad \frac{d\mu(x)}{dx} = V(x)^{-\beta},$$

where $V$ is a positive convex function on $\mathbb{R}^n$, defined on some open convex set in $\mathbb{R}^n$, and $\beta \geq n$. If $\beta = n$, this formula describes the class of all absolutely continuous convex measures on $\mathbb{R}^n$. More precisely, the measure $\mu$ will be $\kappa$-concave with $\kappa = -1/(\beta - n)$, which may or may not be optimal, depending on the choice of $V$.

The Borell description (1.8), or, more precisely its sufficiency part, represents a particular case of the following lemma that describes the Borell–Brascamp–Lieb dimensional extension of the Prékopa–Leindler theorem [9, 10] (cf. [11, 13],...). This result is the one which is suited to the family of convex measures considered here.

LEMMA 1.1. *Let $\kappa \leq 1/n$ and $\kappa_n = \kappa/(1 - n\kappa)$. Given $0 < t < 1$, let $u, v, w$ be nonnegative measurable functions, defined on some open convex set $\Omega \subset \mathbb{R}^n$ and satisfying*

$$(1.10) \qquad w(tx + (1-t)y) \geq [tu(x)^{\kappa_n} + (1-t)v(y)^{\kappa_n}]^{1/\kappa_n}$$

*for all $x, y \in \Omega$ such that $u(x) > 0$ and $v(y) > 0$. Then*

$$(1.11) \qquad \int_\Omega w(z)\,dz \geq \left[t\left(\int_\Omega u(x)\,dx\right)^\kappa + (1-t)\left(\int_\Omega v(y)\,dy\right)^\kappa\right]^{1/\kappa}.$$

This result will be the key tool to the following generalization of the Brascamp–Lieb theorem (1.5). Namely, let $\mu$ be a probability measure of the type (1.9). Under the additional assumption of twice differentiability on $V$, if $\beta \geq n$, for any bounded smooth function $g$ with mean zero,

$$(1.12) \qquad \mathrm{Var}_\mu(g) \leq \frac{1}{\beta + 1} \int \frac{\langle V''^{-1} \nabla G, \nabla G \rangle}{V}\,d\mu,$$

where $G = gV$. Replacing here $V$ with $c_\beta^{1/\beta}(1 + \frac{1}{\beta}W)$, where $c_\beta$ is a normalizing constant, and letting $\beta \to +\infty$, we indeed arrive at (1.5). On the other hand, (1.12) will be shown to easily imply the weighted Poincaré-type inequality (1.2) for the Cauchy distributions.

The extension (1.12) of the Brascamp–Lieb theorem is proved in Section 2 on the basis of Lemma 1.1. It is applied next to derive a quantitative refinement of (1.5) within the class of log-concave measures. In Section 3, the extended Brascamp–Lieb theorem is applied to the Cauchy family to establish in particular (1.2) and (1.3). The family of Cauchy distributions $\nu_\beta$, $\beta \geq (n+1)/2$, also shares another functional inequality of interest, namely the following weighted logarithmic Sobolev inequality:

$$\mathrm{Ent}_{\nu_\beta}(g^2) \leq \frac{1}{\beta - 1} \int |\nabla g(x)|^2 (1 + |x|^2)^2\,d\nu_\beta(x)$$

with thus (and naturally) a worse weight function in comparison with (1.2). Here $\mathrm{Ent}_\mu(g^2) = \int g^2 \log g^2\,d\mu - \int g^2\,d\mu \log \int g^2\,d\mu$ denotes the entropy of $g^2$ under $\mu$. This weighted logarithmic Sobolev inequality is established by a different procedure, as a consequence actually of the Gaussian case.

Abstract weighted Poincaré and logarithmic Sobolev inequalities are connected with the problem of large deviations of Lipschitz functions and measure concentration. This aspect is studied in Section 4 and illustrated on the Cauchy examples. As announced, the tails of Lipschitz functions with respect to $\nu_\beta$ admit in turn Gaussian, exponential, and polynomial decays as the interval tends to infinity.

In the last section, we consider weighted isoperimetric inequalities and show that, somehow surprisingly, all concave probability measures described by (1.9) share the weighted Poincaré-type inequality (1.1) with, up to a constant, the universal weight $w(x)^2 = 1 + |x|^2$. The result is established through a stronger Cheeger-type isoperimetric inequality. This universality property illustrates the importance of this weight and of the family of Cauchy distributions in this investigation.

After this work was completed, we became aware of the related references [2, 3], where some similar inequalities are proved and analyzed by different methods and with different purposes.

**2. Brascamp–Lieb-type inequality.** Consider throughout this section an absolutely continuous probability measure $\mu$, concentrated on an open, convex set $\Omega \subset \mathbb{R}^n$, with density

$$p(x) = V(x)^{-\beta}, \qquad x \in \Omega,$$

where $\beta > n$ and $V$ is a positive, twice continuously differentiable convex function on $\Omega$ with positive second derivative $V''$ in the matrix sense with inverse $V''^{-1}$. We associate to $V$ a distance-like or cost function

$$(2.1) \qquad d_V(x, y) = V(y) - V(x) - \langle V'(x), y - x \rangle, \qquad x, y \in \Omega.$$

Note that $d_V(x, y)$ is nonnegative and is only vanishing when $x = y$. For example, $d_V(x, y) = |x - y|^2$ for $V(x) = 1 + |x|^2$.

As a first step toward (1.12) we prove the following result:

THEOREM 2.1. *Let $\beta > n$, and let $f, g$ be measurable functions on $\Omega$ satisfying, for all $x, y \in \Omega$,*

$$(2.2) \qquad f(x)V(x) \leq g(y)V(y) + d_V(x, y).$$

*If $f$ is $\mu$-integrable and $g \geq -1$ on $\Omega$, then*

$$(2.3) \qquad 1 + \frac{\beta}{\beta - n} \int f \, d\mu \leq \left[ \int (1 + g)^{-\beta} \, d\mu \right]^{-1/(\beta - n)}.$$

The statement remains to hold for arbitrary positive convex differentiable functions $V$ on $\Omega$ with $\int V^{-\beta}\,dx = 1$. One particular case is worth mentioning separately.

COROLLARY 2.2. *Let $\mu$ be a probability measure on $\Omega$ with density $e^{-W}$, where $W$ is a differentiable convex function on $\Omega$. For all measurable functions $f, g$ on $\Omega$ satisfying*

(2.4) $$f(x) \leq g(y) + d_W(x,y)$$

*for all $x, y \in \Omega$ and such that $g$ is $\mu$-integrable,*

(2.5) $$\int e^f \, d\mu \leq e^{\int g \, d\mu}.$$

To deduce the corollary, apply Theorem 2.1 to the functions $f_\beta = -g/\beta$, $g_\beta = -f/\beta$ with $(f,g)$ satisfying (2.4), and to the densities $c_\beta(1 + \frac{1}{\beta}W)^{-\beta}$. That is, with $V = c_\beta^{1/\beta}(1 + \frac{1}{\beta}W)$, $d_V = c_\beta^{1/\beta} d_W/\beta$ in (2.1) and (2.2). Then, under unessential technical assumptions on $f$, $g$ and $W$, the inequality (2.2) will be fulfilled as $\beta \to +\infty$, and in the limit (2.3) will turn into (2.5).

In the Gaussian case, inequality (2.5) was obtained by Tsirelson [25], with convex $g$, who developed Chevet's concept of Gaussian mixed volumes, and by Maurey [20], who proposed a different approach using the infimum-convolution operator and the Prékopa–Leindler theorem. As shown in [6], infimum-convolution estimates for the quadratic cost $d_W$ are connected with logarithmic Sobolev inequalities. Maurey's argument was further developed and extended to the class of log-concave probability measures in [7], where Corollary 2.2 was established for "potentials" $W$ whose associated distance $d_W$ is controlled in terms of a norm on $\mathbb{R}^n$.

As announced, to reach Theorem 2.1, we make use of the fundamental Lemma 1.1 (with $\kappa < 0$).

PROOF OF THEOREM 2.1. One may assume that $f$ is bounded from above and $\inf_\Omega g > -1$. For $t \in (0,1)$, $s = 1 - t$, define $T_s(x,y) = \frac{1}{ts}[tV(x) + sV(y) - V(tx + sy)]$ so that

(2.6) $$\lim_{s \to 0} T_s(x,y) = d_V(x,y).$$

The convergence is uniform on compact subsets of $\Omega \times \Omega$. Take an arbitrary open convex set $\Omega_0$ with compact closure in $\Omega$. Given $t \in (0,1)$ and $\varepsilon > 0$, apply Lemma 1.1 with $\kappa = -1/(\beta - n)$, $\kappa_n = \kappa/(1 - n\kappa) = -1/\beta$, to the functions on $\Omega_0$,

$$u(x) = (1 - sf_\varepsilon(x))^{1/\kappa_n} p(x)/\mu(\Omega_0),$$
$$v(y) = (1 + tg(y))^{1/\kappa_n} p(y)/\mu(\Omega_0),$$
$$w(z) = p(z)/\mu(\Omega_0),$$

where $f_\varepsilon = f - \varepsilon/V$. Note that $u$ is well-defined and positive for sufficiently small $s > 0$. Now, condition (1.10) in Lemma 1.1 reads as

$$f(x)V(x) \leq g(y)V(y) + T_s(x,y) + \varepsilon, \qquad x, y \in \Omega_0.$$

Due to the hypothesis (2.2) and the property (2.6), this condition is fulfilled for all sufficiently small $s > 0$ uniformly on $\Omega_0$. Hence, we obtain (1.11) of Lemma 1.1, that is,

$$(2.7) \qquad 1 \leq t\left[\int (1 - sf_\varepsilon)^{1/\kappa_n}\, d\nu_0\right]^\kappa + s\left[\int (1 + tg)^{1/\kappa_n}\, d\nu_0\right]^\kappa,$$

where $\nu_0$ denotes the normalized restriction of $\mu$ to the set $\Omega_0$. By Taylor's expansion, (2.7) yields in the limit as $s \to 0$,

$$1 + \frac{\kappa}{\kappa_n}\int f_\varepsilon\, d\nu_0 \leq \left[\int (1+g)^{1/\kappa_n}\, d\nu_0\right]^\kappa.$$

It remains to let first $\varepsilon \to 0$ and then $\Omega_0 \uparrow \Omega$ to get the desired inequality (2.3). The proof of Theorem 2.1 is complete. $\square$

On the basis of Theorem 2.1, we establish our main extension of the Brascamp–Lieb theorem.

THEOREM 2.3. *Let $\beta > n$. For any smooth, $\mu$-square integrable function $g$ on $\Omega$, and $G = Vg$,*

$$(2.8) \qquad (\beta+1)\operatorname{Var}_\mu(g) \leq \int \frac{\langle V''^{-1}\nabla G, \nabla G\rangle}{V}\, d\mu + \frac{n}{\beta - n}\left(\int g\, d\mu\right)^2.$$

If we assume that $g$ has mean zero, so that the last term in (2.8) is vanishing, we arrive at the announced bound (1.12) from the Introduction. At this step, the condition $\beta > n$ may be relaxed to $\beta \geq n$ by continuity.

PROOF OF THEOREM 2.3. We may assume that $\|V''\|$ and $\|V''^{-1}\|$ are bounded and uniformly continuous on $\Omega$. Assume also $g$, $|\nabla G|$ and $\|G''\|$ are bounded [otherwise, restrict all inequalities to open, convex sets $\Omega_0$ with compact closure in $\Omega$, and then approximate the latter with $\Omega_0$'s in the resulting inequality (2.8)]. Given $\varepsilon > 0$, define

$$F_\varepsilon(x) = \inf_{y \in \Omega}\left[G(y) + \frac{1}{\varepsilon}d_V(x,y)\right], \qquad x \in \Omega$$

with $d_V$ introduced in (2.1), and set $f_\varepsilon = F_\varepsilon/V$. That is, $\varepsilon f_\varepsilon$ is optimal for the function $\varepsilon g$ in the "infimum-convolution" inequality (2.2). Note that $F_\varepsilon$ is upper-semicontinuous as infimum of a family of continuous functions. By Theorem 2.1 applied to the couple $(\varepsilon f_\varepsilon, \varepsilon g)$ with sufficiently small $\varepsilon$ (when $\varepsilon g \geq -1$),

$$(2.9) \qquad 1 + \varepsilon \frac{\kappa}{\kappa_n} \int f_\varepsilon \, d\mu \leq \left[ \int (1+\varepsilon g)^{1/\kappa_n} \, d\mu \right]^\kappa,$$

where $\kappa = -1/(\beta - n)$ and $\kappa_n = -1/\beta$.

We claim that inequality (2.9), as $\varepsilon \to 0$, yields (2.8). Note, by the assumption on $V$, we have $d(x, x+\varepsilon h) \geq c|h|^2 \varepsilon^2$ with a constant $c > 0$, depending $V$, only. Hence, the infimum in the definition of $F_\varepsilon(x)$ may be restricted to the points $y = x + \varepsilon h$ with $|h| < r$, where $r$ depends on $V$ and $G$. Moreover, by Taylor's expansion

$$d_V(x, x+\varepsilon h) = \frac{\varepsilon^2}{2} \langle V''(x)h, h \rangle + |h|^2 o(\varepsilon^2),$$

where the constant in $o(\varepsilon^2)$ is numerical, that is, it may be chosen to be independent of $(x, h)$ by the uniform continuity of $V''$. Hence, together with the Taylor expansion for $G(x + \varepsilon h)$, we get that

$$\begin{aligned}
F_\varepsilon(x) &= \inf_{h\,:\,x+\varepsilon h \in \Omega} \left[ G(x+\varepsilon h) + \frac{1}{\varepsilon} d_V(x, x+\varepsilon h) \right] \\
&= \inf_{h\,:\,x+\varepsilon h \in \Omega} \left[ G(x) + \varepsilon \langle \nabla G(x), h \rangle + \frac{\varepsilon}{2} \langle V''(x)h, h \rangle + |h|^2 o(\varepsilon) \right] \\
&\geq \inf_{h \in \mathbb{R}^n} \left[ G(x) + \varepsilon \langle \nabla G(x), h \rangle + \frac{\varepsilon}{2} \langle V''(x)h, h \rangle + |h|^2 o(\varepsilon) \right] \\
&= G_\varepsilon(x) + o(\varepsilon),
\end{aligned}$$

where $G_\varepsilon(x) = G(x) - \frac{\varepsilon}{2} \langle V''^{-1} \nabla G(x), \nabla G(x) \rangle$. Hence,

$$(2.10) \qquad f_\varepsilon(x) = g(x) - \frac{\varepsilon}{2} \frac{\langle V''^{-1} \nabla G(x), \nabla G(x) \rangle}{V(x)} + o(\varepsilon).$$

Now, let us turn to the right-hand side of (2.9). By Taylor's expansion, using the expectation sign $\mathbb{E}$ for integrals over the measure $\mu$,

$$[\mathbb{E}(1+\varepsilon g)^{1/\kappa_n}]^\kappa = 1 + \frac{\varepsilon \kappa}{\kappa_n} \mathbb{E}g + \frac{\varepsilon^2}{2} \frac{\kappa}{\kappa_n} \left( \frac{1}{\kappa_n} - 1 \right) \mathbb{E}g^2 + \frac{\kappa(\kappa-1)}{2} \left( \frac{\varepsilon}{\kappa_n} \mathbb{E}g \right)^2 + o(\varepsilon^2).$$

Subtracting 1 from both sides of (2.9) and dividing by $\varepsilon^2 \frac{\kappa}{\kappa_n}$, we arrive together with the latter at

$$\frac{\mathbb{E}f_\varepsilon - \mathbb{E}g}{\varepsilon} \leq \frac{1}{2} \left[ \left( \frac{1}{\kappa_n} - 1 \right) \mathbb{E}g^2 + \frac{\kappa - 1}{\kappa_n} (\mathbb{E}g)^2 \right] + o(1).$$

Finally, by (2.10), comparing the coefficients in front of $\varepsilon$, we get

$$\left( 1 - \frac{1}{\kappa_n} \right) \mathbb{E}g^2 \leq \mathbb{E} \frac{\langle V''^{-1} \nabla G, \nabla G \rangle}{V} + \frac{\kappa - 1}{\kappa_n} (\mathbb{E}g)^2,$$

which is an equivalent form of (2.8). Theorem 2.3 has thus been proved. □

We suggest an alternate formulation of Theorem 2.3. Given $V$, a positive twice continuously differentiable convex function on some open convex set $\Omega$ in $\mathbb{R}^n$, let $\mu_\beta$ be the probability measure with density $p(x) = Z_\beta^{-1} V(x)^{-\beta}$, $\beta > n$, where it will be useful here to specify the normalization $Z_\beta$. If we rewrite Theorem 2.3 in terms of $G$ rather than $g$, we see that whenever $\int G \, d\mu_{\beta+1} = 0$, then

$$(1+\beta) \int \frac{G^2}{V^2} \, d\mu_\beta \leq \int \frac{\langle V''^{-1} \nabla G, \nabla G \rangle}{V} \, d\mu_\beta.$$

By definition of $\mu_\beta$, and changing $\beta$ into $\beta - 1 \geq n$,

$$(2.11) \qquad \int G^2 \, d\mu_{\beta+1} \leq \frac{1}{\beta} \cdot \frac{Z_\beta}{Z_{\beta+1}} \int \langle V''^{-1} \nabla G, \nabla G \rangle \, d\mu_\beta$$

provided that $\int G \, d\mu_\beta = 0$. This provides an alternate extension of the Brascamp–Lieb inequality (1.5). In particular, if $V'' \geq c \operatorname{Id}$ for some $c > 0$ in the sense of symmetric matrices, and $\beta \geq n+1$,

$$(2.12) \qquad \int G^2 \, d\mu_{\beta+1} \leq \frac{1}{c\beta} \cdot \frac{Z_\beta}{Z_{\beta+1}} \int |\nabla G|^2 \, d\mu_\beta$$

for all smooth $G$ with $\int G \, d\mu_\beta = 0$, as an analogue of (1.6).

To conclude this section, we briefly mention an application of Theorem 2.3 to some improved bound for logarithmically concave measures. Thus, let $\mu$ be a probability measure on an open convex set $\Omega$ in $\mathbb{R}^n$ with density $e^{-W}$, where $W$ is a function on $\Omega$ of class $C^2$ with positive second derivative. In (2.8) we may take $V = e^{W/\beta}$, which has the first two derivatives $V' = \frac{V}{\beta} W'$ and $V'' = \frac{V}{\beta}(W'' + \frac{1}{\beta} W' \otimes W')$. Here, for vectors $v \in \mathbb{R}^n$, we use the tensor product notation $v \otimes v$ to denote the $n \times n$ matrix with entries $v_i v_j$.

Let $g$ be smooth and such that $\int g \, d\mu = 0$. To estimate the integrand on the right of (2.8), write $(gV)' = gV' + Vg'$ (where we simplify notations from $\nabla g$ to $g'$). Apply the elementary bound

$$(2.13) \quad \langle A(u+v), (u+v) \rangle \leq r \langle Au, u \rangle + \frac{r}{r-1} \langle Av, v \rangle, \qquad u, v \in \mathbb{R}^n, r > 1,$$

where $A$ is an arbitrary positively definite, symmetric matrix, to get that

$$(2.14) \quad \langle V''^{-1}(gV)', (gV)' \rangle \leq r \langle V''^{-1} V', V' \rangle g^2 + \frac{r}{r-1} \langle V''^{-1} g', g' \rangle V^2.$$

Now, by the convexity of $W$,

$$\frac{\langle V''^{-1} V', V' \rangle}{V} = \langle (\beta W'' + W' \otimes W')^{-1} W', W' \rangle \leq 1.$$

Hence, applying (2.14) in (2.8) and introducing the family of positively definite matrices

(2.15) $$R_{W,\beta} = W'' + \frac{1}{\beta} W' \otimes W',$$

we obtain that

$$(\beta + 1 - r) \int g^2 \, d\mu \leq \frac{r\beta}{r-1} \int \langle R_{W,\beta}^{-1} \nabla g, \nabla g \rangle \, d\mu.$$

It remains to optimize over all $r > 1$ to conclude with the following statement:

THEOREM 2.4. *For any smooth, $\mu$-square integrable function $g$ on $\Omega$, and for any $\beta \geq n$,*

(2.16) $$\mathrm{Var}_\mu(g) \leq C_\beta \int \langle R_{W,\beta}^{-1} \nabla g, \nabla g \rangle \, d\mu,$$

*where $C_\beta = (1 + \sqrt{\beta+1})^2/\beta$.*

Note $1 < C_\beta < 6$ and $C_\beta \to 1$ as $\beta$ grows to infinity. Since $R_{W,\beta}^{-1} \leq W''^{-1}$ in the matrix sense, the Brascamp–Lieb inequality (1.5) may thus be viewed as a limiting case of (2.16). On the other hand, finite values of $\beta$ may give an improvement in terms of the decay of the weight function. In particular, in dimension one with $\beta = 1$ we always have

(2.17) $$\mathrm{Var}_\mu(g) \leq 6 \int \frac{g'(x)^2}{W''(x) + W'(x)^2} \, d\mu(x).$$

For example, when $\mu$ has density $p(x) = \lambda e^{-\lambda x}$ on $\Omega = (0, +\infty)$ with a positive parameter $\lambda$, we arrive at the usual Poincaré-type inequality $\mathrm{Var}_\mu(g) \leq \frac{6}{\lambda^2} \int g'^2 \, d\mu$, which cannot be obtained on the basis of (1.5). For the Gaussian measure $\mu = \gamma_1$, (2.17) gives a weighted Poincaré-type inequality

$$\mathrm{Var}_{\gamma_1}(g) \leq 6 \int \frac{g'(x)^2}{1 + x^2} \, d\gamma_1(x)$$

with an asymptotically sharp weight function. Note, however, for the $n$-dimensional Gaussian measure $\mu = \gamma_n$ we only have from (2.16) with $\beta = n$,

(2.18) $$\mathrm{Var}_{\gamma_n}(g) \leq 6 \int \left[ |\nabla g(x)|^2 - \frac{\langle \nabla g(x), x \rangle^2}{n + |x|^2} \right] d\gamma_n(x).$$

Being restricted to radial functions $g(x) = g(|x|)$, the latter yields a weighted Poincaré-type inequality for the family of $\chi_n$-distributions.

**3. The generalized Cauchy distribution.** In this section, we specialize Theorem 2.3 to the case of the generalized Cauchy distributions $\nu_\beta$ and develop for this specific family an analytic and geometric investigation similar to the one for the classical Gaussian model. Recall the generalized Cauchy distribution $\nu_\beta$ has the density

$$p(x) = \frac{1}{Z}(1+|x|^2)^{-\beta}, \qquad x \in \mathbb{R}^n, \ \beta > \frac{n}{2}.$$

Here, $Z = n\omega_n \Gamma(\frac{n}{2})\Gamma(\beta - \frac{n}{2})/2\Gamma(\beta)$ is a normalizing factor (where $\omega_n$ denotes the volume of the Euclidean ball in $\mathbb{R}^n$ of radius one). But it has no influence for the first integral in (2.8), so one may take $V(x) = 1 + |x|^2$. In this case $V''^{-1} = \frac{1}{2}\operatorname{Id}$, and if $g$ has $\nu_\beta$-mean zero, Theorem 2.3 with $\beta \geq n$ tells us that

(3.1) $$(\beta+1)\int g^2\, d\nu_\beta \leq \int \frac{|\nabla G(x)|^2}{2(1+|x|^2)}\, d\nu_\beta(x),$$

where $G(x) = g(x)(1+|x|^2)$. Evidently, $|\nabla G(x)| \leq 2|g(x)||x| + |\nabla g(x)|(1+|x|^2)$, and applying (2.13) with parameter $r > 1$,

$$|\nabla G(x)|^2 \leq 4rg(x)^2|x|^2 + \frac{r}{r-1}|\nabla g(x)|^2(1+|x|^2)^2.$$

Using this estimate in (3.1) we obtain a family of the weighted Poincaré-type inequalities

$$\int g^2\, d\nu_\beta \leq \frac{r}{2(r-1)((\beta+1)-2r)} \int |\nabla g(x)|^2 (1+|x|^2)\, d\nu_\beta(x).$$

It is easy to check that the constant on the right is minimized for $r = \sqrt{(\beta+1)/2}$ and then we arrive at:

THEOREM 3.1. *The generalized Cauchy distribution $\nu_\beta$ with $\beta \geq n$ satisfies the weighted Poincaré-type inequality*

(3.2) $$\operatorname{Var}_{\nu_\beta}(g) \leq \frac{C_\beta}{2(\beta-1)} \int |\nabla g(x)|^2(1+|x|^2)\, d\nu_\beta(x)$$

*for all bounded smooth functions $g$ on $\mathbb{R}^n$ with $C_\beta = (\sqrt{1+\frac{2}{\beta-1}} + \sqrt{\frac{2}{\beta-1}})^2$.*

Note that $C_\beta > 1$ and $C_\beta \to 1$ as $\beta \to +\infty$. So, after the linear change of the variable $y = \sqrt{2\beta}x$, in the limit we will be led in (3.2) to the Gaussian Poincaré-type inequality with optimal constant.

To get more information on how optimal the constant in (3.2) is, one can just test the weighted Poincaré-type inequality

(3.3) $$\operatorname{Var}_{\nu_\beta}(g) \leq C(\beta,n) \int |\nabla g(x)|^2(1+|x|^2)\, d\nu_\beta(x)$$

on simple functions. Take, for example, $g(x) = 1/(1+|x|^2)$. For every $m = 1, 2, \ldots$,

$$I_m \equiv \int (1+|x|^2)^{-m} d\nu_\beta(x) = \frac{(\beta - n/2)(\beta - n/2 + 1) \cdots (\beta - n/2 + (m-1))}{\beta(\beta+1) \cdots (\beta + (m-1))}.$$

Now, $\mathrm{Var}_{\nu_\beta}(g) = I_2 - I_1^2 = (\beta - \frac{n}{2})\frac{n}{2}/\beta^2(\beta+1)$. Since $|\nabla g(x)|^2 = 4|x|^2/(1+|x|^2)^4$, the integral in (3.3) equals $4(I_2 - I_3) = 4(\beta - \frac{n}{2})(\beta - \frac{n}{2} + 1)\frac{n}{2}/\beta(\beta+1)(\beta+2)$. Comparing both sides, we conclude that

$$C(\beta, n) \geq \frac{\beta + 2}{4\beta(\beta - n/2 + 1)} \geq \frac{1}{4\beta}.$$

Thus, the constant in (3.2) is optimal within universal factors at least in the region $\beta \geq \max\{n, 2\}$. As for the region $n/2 < \beta < n$, the optimal value of $C(\beta, n)$ essentially depends on the dimension $n$ (see below).

Let us mention that in dimension one this constant may be analyzed directly by reducing the weighted Poincaré-type inequality (3.3) to the Hardy-type inequality

$$\int_0^{+\infty} g(x)^2 p(x)\, dx \leq C \int_0^{+\infty} g'(x)^2 q(x)\, dx$$

with specific weights $p(x) = (1+|x|^2)^{-\beta}$, $q(x) = (1+|x|^2)^{-(\beta-1)}$ (where $g$ is smooth on $[0, +\infty)$ with $g(0) = 0$). In general (cf. [21]), one has $B \leq C \leq 4B$, where $B = \sup_{x>0}[\int_0^x dt/q(t) \int_x^{+\infty} p(t)\, dt]$. By elementary calculus, this quantity may be bounded in the weighted Cauchy case by $1/\max\{2(\beta-1), 1\}$. Hence,

$$C(\beta, 1) \leq \frac{4}{\max\{2(\beta-1), 1\}} \leq \frac{6}{\beta}$$

in the whole range $\beta > 1$. Together with (3.2) for the case $n \geq 2$, we obtain that $C(\beta, n) \leq C/\beta$ for all $\beta \geq n$ ($\beta > 1$) with some numerical constant $C$.

The following corollary is concerned with the reversed form (1.3).

COROLLARY 3.2. *If $\beta \geq n+1$, for any smooth bounded function $g$ on $\mathbb{R}^n$,*

(3.4) $$\inf_{c \in \mathbb{R}} \int \frac{|g(x) - c|^2}{1 + |x|^2}\, d\nu_\beta(x) \leq \frac{1}{2\beta} \int |\nabla g|^2\, d\nu_\beta.$$

The proof of Corollary 3.2 is an immediate consequence of (2.12) (or directly Theorem 2.3) with $c = 2$. The remaining range $n/2 < \beta \leq n+1$ in Corollary 3.2 may be treated similarly by choosing a different "potential," $V_\alpha(x) = (1+|x|^2)^\alpha$, with $1/2 < \alpha \leq 1$. At every point $x \in \mathbb{R}$, $x \neq 0$, it has a positive Hessian $V''_\alpha(x) = 2\alpha(1+|x|^2)^{\alpha-1}(\mathrm{Id} - \lambda e \otimes e)$, where $e = x/|x|$, $\lambda = 2(1-\alpha)|x|^2/(1+|x|^2)$, with the inverse matrix

$$V''_\alpha(x)^{-1} = \frac{1}{2\alpha(1+|x|^2)^{\alpha-1}}\left(\mathrm{Id} + \frac{\lambda}{1-\lambda}e \otimes e\right).$$

It follows that $\|V''_\alpha(x)^{-1}\| \leq 1/[2\alpha(2\alpha-1)(1+|x|^2)^{\alpha-1}]$.

Now, given $\beta > n/2$, write $\beta = 2\alpha + \beta' - 1$ with $\beta' = \alpha\gamma$, $\gamma \geq n$. Applying Theorem 2.3 to $\nu_{\beta'}$ with density written as $\frac{1}{Z}V_\alpha^{-\gamma}$ (i.e., with $\gamma$ in place of $\beta$), we obtain similarly to (3.1) that

$$(\gamma + 1)\int g^2\, d\nu_{\beta'} \leq \frac{1}{2\alpha(2\alpha-1)}\int \frac{|\nabla G(x)|^2}{(1+|x|^2)^{2\alpha-1}}\, d\nu_{\beta'}(x),$$

where $g$ is bounded, smooth, with $\nu_{\beta'}$-mean zero, and $G = V_\alpha g$. Equivalently, in terms of $\nu_\beta$ and $G$, we have that $(\gamma+1)\int \frac{|G(x)|^2}{1+|x|^2}\, d\nu_\beta(x) \leq \frac{1}{2\alpha(2\alpha-1)}\int |\nabla G|^2\, d\nu_\beta$. Changing $G$ into $g$, we deduce that

$$\inf_{c\in\mathbb{R}}\int \frac{|g(x) - c|^2}{1+|x|^2}\, d\nu_\beta(x) \leq \frac{1}{2\alpha(2\alpha-1)(\gamma+1)}\int |\nabla g|^2\, d\nu_\beta.$$

If $n/2 < \beta \leq n+1$, one may just choose $\alpha = (\beta+1)/(n+2)$, which leads to the analogue of the reversed weighted Poincaré-type inequality (3.4),

$$(3.5) \qquad \inf_{c\in\mathbb{R}}\int \frac{|g(x) - c|^2}{1+|x|^2}\, d\nu_\beta(x) \leq C\int |\nabla g|^2\, d\nu_\beta$$

with constant

$$C = \frac{(n+2)^2}{2(n+1)(\beta+1)}\frac{1}{2\beta - n} \leq \frac{1}{\beta - n/2}.$$

A similar approach may also be used to involve the values $n/2 < \beta \leq n$ in Theorem 3.1, but we leave this to the reader as an exercise. Instead, let us note that the reversed form, such as (3.4) and (3.5), can be deduced from the weighted Poincaré-type inequalities, such as (3.2). Namely, in place of the variance in (3.1), one may consider other similar quantities. For example, due to the elementary general bound $\mathrm{Var}(g) \geq \frac{1}{3}\mathbb{E}_\mu|g - m|^2$, where $m$ is a median of $g$ with respect to $\mu$, (3.2) yields

$$\int |g - m|^2\, d\nu_\beta \leq \frac{3C_\beta}{2(\beta - 1)}\int |\nabla g(x)|^2(1+|x|^2)\, d\nu_\beta(x).$$

With the help of this inequality, Corollary 3.2 immediately follows for sufficiently large $\beta$ (say, $\beta \geq 7$) and with an additional numerical factor in view of the following general proposition of possible independent interest.

PROPOSITION 3.3. *Assume a probability measure $\mu$ on $\mathbb{R}^n$ satisfies the weighted Poincaré-type inequality*

$$\text{(3.6)} \qquad \int |g - m|^2 \, d\mu \leq \int |\nabla g(x)|^2 (a + b|x|^2) \, d\mu(x)$$

*with some constants $a > 0$ and $b \in [0, 1)$. Then for any smooth $g$ on $\mathbb{R}^n$,*

$$\text{(3.7)} \qquad \inf_{c \in \mathbb{R}} \int \frac{|g(x) - c|^2}{a + b|x|^2} \, d\mu(x) \leq \frac{1}{(1 - \sqrt{b})^2} \int |\nabla g|^2 \, d\mu.$$

PROOF. Indeed, we may restrict ourselves to nonnegative $g$ with median zero. Fix such a function and consider $f(x) = g(x)/\sqrt{a + b|x|^2}$. By the one-dimensional variant of (2.13), for every $r > 1$,

$$|\nabla f(x)|^2 \leq \frac{r}{r-1} \frac{|\nabla g(x)|^2}{a + b|x|^2} + rb \frac{g(x)^2}{(a + b|x|^2)^2}.$$

This estimate may be used in (3.6) with $f$ in place of $g$ to get that, whenever $rb < 1$,

$$\int \frac{g(x)^2}{a + b|x|^2} \, d\mu(x) \leq \frac{r}{(r-1)(1 - rb)} \int |\nabla g(x)|^2 \, d\mu(x).$$

The optimal choice $r = 1/\sqrt{b}$ then leads to the conclusion (3.7). □

Returning to the reversed form (3.4) and (3.5), let us finally note that the weight function on the left-hand side may be slightly improved: one of the results by Denzler and McCann in [12] states that, for $n \geq 2$,

$$\inf_{c \in \mathbb{R}} \int \frac{|g(x) - c|^2}{|x|^2} \, d\nu_\beta(x) \leq C(\beta, n) \int |\nabla g|^2 \, d\mu.$$

As announced in the introductory section, in analogy with the Gaussian measure, the Cauchy measures $\nu_\beta$ admit a weighted logarithmic Sobolev inequality, but with a different weight function in comparison with (3.1).

THEOREM 3.4. *If $\beta \geq (n+1)/2$, $\beta > 1$, for any smooth bounded $g$ on $\mathbb{R}^n$,*

$$\text{(3.8)} \qquad \operatorname{Ent}_{\nu_\beta}(g^2) \leq \frac{1}{\beta - 1} \int |\nabla g(x)|^2 (1 + |x|^2)^2 \, d\nu_\beta(x).$$

After linear rescaling of the coordinates, the inequality (3.5) with growing $\beta$ yields the Gross logarithmic Sobolev inequality for the Gaussian measure $\gamma_n$ on $\mathbb{R}^n$,

$$\text{(3.9)} \qquad \operatorname{Ent}_{\gamma_n}(g^2) \leq 2 \int |\nabla g|^2 \, d\gamma_n.$$

However, the proof of Theorem 3.4 may be given on the basis of this limiting case itself.

PROOF OF THEOREM 3.4. Write $\beta = (n+d)/2$ with $d \geq 1$. The measure $\nu_\beta$ may be characterized as the distribution of the random vector $X = Y/\xi$, where $Y$ is a random vector in $\mathbb{R}^n$ with the standard Gaussian distribution, and $\xi > 0$ is a random variable independent of $Y$ and having the $\chi_d$-distribution with $d$ degrees of freedom. That is, $\xi$ has density

$$\chi_d(r) = \frac{1}{2^{d/2-1}\Gamma(d/2)} r^{d-1} e^{-r^2/2}, \qquad r > 0.$$

In other words, $\nu_\beta$ represents the image of the product probability measure $P = \gamma_n \otimes \chi_d$ on $\mathbb{R}^n \times (0, +\infty)$ under the map $(y, r) \to y/r$ [where, with some abuse, we denote by $\chi_d$ the probability measure with density $\chi_d(r)$].

Since $d \geq 1$, the density $\chi_d$ is log-concave with respect to $\gamma_1$. So, the Bakry–Emery criterion may be applied in dimension one to get the one-dimensional analogue of (3.9), that is, $\mathrm{Ent}_{\chi_d}(u^2) \leq 2\int |u'|^2 \, d\chi_d$ (cf. [17]). Therefore, by the basic product property of logarithmic Sobolev inequalities, the measure $P$ also satisfies the logarithmic Sobolev inequality

$$\mathrm{Ent}_P(f^2) \leq 2\int |\nabla f|^2 \, dP$$

in the class of all smooth functions $f$ on $\mathbb{R}^n \times (0, +\infty)$. Now, given a smooth $g$ on $\mathbb{R}^n$, apply this inequality to functions of the form $f(y, r) = g(y/r)$. Then we get

(3.10) $$\mathrm{Ent}_{\nu_\beta}(g^2) = \mathrm{Ent}_P(f^2) \leq 2 \int |\nabla f|^2 \, dP$$

with $|\nabla f|^2 = |\nabla_y f|^2 + |\nabla_r f|^2$. Letting $x = y/r$, we have

$$\frac{\partial f(y, r)}{\partial y_i} = \frac{1}{r} \frac{\partial g(x)}{\partial y_i}, \qquad i = 1, \ldots, n, \qquad \frac{\partial f(y, r)}{\partial r} = -\frac{1}{r} \langle \nabla g(x), x \rangle,$$

so that $|\nabla f(y,r)|^2 \leq \frac{1+|x|^2}{r^2} |\nabla g(x)|^2$. We apply this bound to the right-hand side of (3.10). Namely, changing the variable $y = rx$ and then $r = t/\sqrt{1+|x|^2}$, we get in terms of $\Psi(x) = (1+|x|^2)|\nabla g(x)|^2$ that

$$\int |\nabla f|^2 \, dP \leq \int \frac{1}{r^2} \Psi(x) \, dP(y, r)$$

$$= \frac{1}{2^{d/2-1}\Gamma(d/2)(2\pi)^{n/2}}$$

$$\times \int_{\mathbb{R}^n} \frac{\Psi(x)}{(1+|x|^2)^{(n+d)/2-1}} \, dx \int_0^{+\infty} t^{n+d-3} e^{-t^2/2} \, dt.$$

The last expression can be recognized as

$$\int (1+|x|^2)\Psi(x)\,d\nu_\beta(x) \cdot \frac{\int_0^{+\infty} t^{n+d-3} e^{-t^2/2}\,dt}{\int_0^{+\infty} t^{n+d-1} e^{-t^2/2}\,dt}, \tag{3.11}$$

since repeating the previous arguments we also have that

$$1 = \int dP(y,r)$$
$$= \frac{1}{2^{d/2-1}\Gamma(d/2)(2\pi)^{n/2}} \int_{\mathbb{R}^n} \frac{dx}{(1+|x|^2)^{(n+d)/2}} \int_0^{+\infty} t^{n+d-1} e^{-t^2/2}\,dt.$$

But $\int_0^{+\infty} t^{d-1} e^{-t^2/2}\,dt = 2^{d/2-1}\Gamma(d/2)$, so the ratio in (3.11) is equal to

$$\frac{2^{(n+d-2)/2-1}\Gamma((n+d-2)/2)}{2^{(n+d)/2-1}\Gamma((n+d)/2)} = \frac{1}{2}\frac{1}{(n+d)/2-1} = \frac{1}{2(\beta-1)}.$$

It remains to combine (3.10) and (3.11) to conclude the argument. The proof of Theorem 3.4 is complete. □

**4. Growth of moments and large deviations.** Once we have realized that the generalized Cauchy distributions satisfy weighted Poincaré-type and weighted logarithmic Sobolev inequalities

$$\mathrm{Var}_\mu(g) \leq \int |\nabla g|^2 w^2 \,d\mu, \tag{4.1}$$

$$\mathrm{Ent}_\mu(g) \leq 2 \int |\nabla g|^2 w^2 \,d\mu \tag{4.2}$$

with some specific weight functions $w^2$, it is natural to wonder what kind of information may be deduced from these functional inequalities themselves. In particular, one is typically interested in the moment and large deviation bounds for Lipschitz functions, parts of the concentration of measure phenomenon. As we will see, the Gromov–Milman theorem and the so-called Herbst argument (cf. [14, 16, 17]) may be adapted in a natural way to the case of a general weight $w^2$ to produce probability decays that fit the nature of the Cauchy and more general convex distributions.

Given a measurable function $f$ on $\mathbb{R}^n$, define $L^p$-norms or $p$th moments with respect to $\mu$ by $\|f\|_p^p = \mathbb{E}|f|^p = \int |f|^p\,d\mu$, where $p \geq 1$.

THEOREM 4.1. *Assume $w$ has a finite $p$th moment, $p \geq 2$. Then under (4.1) any Lipschitz function $f$ on $\mathbb{R}^n$ has a finite $p$th moment, and if $\|f\|_{\mathrm{Lip}} \leq 1$, $\int f\,d\mu = 0$,*

$$\|f\|_p \leq \frac{p}{\sqrt{2}} \|w\|_p. \tag{4.3}$$

*Moreover, under (4.2),*

$$\|f\|_p \leq \sqrt{p-1}\|w\|_p. \tag{4.4}$$

PROOF. Let $\|f\|_{\text{Lip}} \leq 1$ and $\int f \, d\mu = 0$. We may assume $f$ is bounded (otherwise apply a simple truncation argument). The inequality (4.1) can be tensorized to yield, for any bounded measurable function $g$ on $\mathbb{R}^n \times \mathbb{R}^n$, which is locally Lipschitz with respect to every variable,

$$\text{Var}_{\mu \otimes \mu}(g) \leq \int [|\nabla_x g(x,y)|^2 w(x)^2 + |\nabla_y g(x,y)|^2 w(y)^2] \, d\mu(x) \, d\mu(y).$$

Since $p \geq 2$, we may apply it to $g(x,y) = \text{sign}(f(x) - f(y))|f(x) - f(y)|^{p/2}$, which gives, due to the Lipschitz property of $f$,

$$\mathbb{E}|f(x) - f(y)|^p \leq \frac{p^2}{4} \mathbb{E}(|f(x) - f(y)|^{p-2}[w(x)^2 + w(y)^2]),$$

where for short we use the expectation sign for integrals with respect to $\mu \otimes \mu$. By Hölder's inequality, the expectation on the right-hand side may be bounded by

$$(\mathbb{E}|f(x) - f(y)|^p)^{(p-2)/p} (\mathbb{E}[w(x)^2 + w(y)^2]^{p/2})^{2/p}$$
$$\leq 2(\mathbb{E}|f(x) - f(y)|^p)^{(p-2)/p} \|w^2\|_{p/2}.$$

Hence, since $\mathbb{E}f = 0$, $\|f\|_p^2 \leq (\mathbb{E}|f(x) - f(y)|^p)^{2/p} \leq \frac{p^2}{2}\|w^2\|_{p/2}$ and (4.3) follows.

As for (4.4), the argument below essentially appears in the paper [1] by Aida and Stroock within the scheme of the usual logarithmic Sobolev inequality. Namely, apply (4.2) to $g = |f|^{p/2}$, so that $|\nabla g| \leq \frac{p}{2}|f|^{p/2-1}$ and

$$\text{Ent}(|f|^p) \leq \frac{p^2}{2} \int |f|^{p-2} w^2 \, d\mu. \tag{4.5}$$

By Hölder's inequality, $\mathbb{E}|f|^{p-2}w^2 \leq (\mathbb{E}|f|^p)^{1-2/p}(\mathbb{E}w^p)^{2/p} = \|f\|_p^{p-2}\|w\|_p^2$. Hence, from (4.5),

$$\frac{d}{dp}\|f\|_p^2 = 2\|f\|_p^2 \frac{\text{Ent}_\mu(|f|^p)}{p^2 \mathbb{E}|f|^p} \leq \|w\|_p^2,$$

and after integration $\|f\|_p^2 - \|f\|_2^2 \leq (p-2)\|w\|_p^2$. Since (4.2) is stronger than (4.1), we also have $\|f\|_2 \leq \|w\|_2 \leq \|w\|_p$. The two bounds imply $\|f\|_p^2 \leq (p-1)\|w\|_p^2$ which is the result. $\square$

In further applications, Theorem 4.1 will be used through the following consequence.

COROLLARY 4.2. *Assume $\|w\|_p \leq C$ for some $p \geq 2$. Under (4.1), for any $f$ on $\mathbb{R}^n$ with $\|f\|_{\mathrm{Lip}} \leq 1$, $\int f \, d\mu = 0$,*

$$\mu(|f| \geq t) \leq \begin{cases} 2e^{-t/Ce}, & \text{if } 0 \leq t \leq t_1, \\ 2\left(\dfrac{Cp}{t}\right)^p, & \text{if } t \geq t_1, \end{cases}$$

*where $t_1 = Cep$. Moreover, under (4.2), in the interval $0 \leq t < t_0$, $t_0 = C\sqrt{ep}$,*

$$\mu(|f| \geq t) \leq 2e^{-t^2/2C^2 e}.$$

Thus, if $C$ is of order 1 and $p$ is large, we still have an exponential bound on a long interval with length proportional to $p$, like in the presence of the usual Poincaré-type inequality, and a Gaussian bound on an interval with length proportional to $\sqrt{p}$, like in the case of the usual logarithmic Sobolev inequality.

PROOF OF COROLLARY 4.2. By Theorem 4.1, if $2 \leq q \leq p$, we have $\|f\|_q \leq Cq$, that is, $\mathbb{E}|f|^q \leq (Cq)^q$, where the expectation is with respect to $\mu$. In the range $0 \leq q \leq 2$, we may use $\|f\|_q \leq \|f\|_2 \leq \|w\|_2 \leq C$, which implies $\mathbb{E}|f|^q \leq 2(Cq)^q$. Indeed, on the positive half-axis $q > 0$, the function $q \to 2q^q$ is minimized at $q = 1/e$, with minimum value $2e^{-1/e} > 1$. Now, by Chebyshev's inequality, for any $0 < q \leq p$ and $t > 0$,

$$(4.6) \qquad \mu(|f| \geq t) \leq \frac{\mathbb{E}|f|^q}{t^q} \leq 2\left(\frac{Cq}{t}\right)^q.$$

Substituing $s = t/C$, write $\mu(|f| \geq t) \leq 2e^{-\varphi(q)}$, where $\varphi(q) = q \log s - q \log q$. This function is concave on $(0, +\infty)$ and attains its maximum $\varphi(q_0) = s/e = t/(Ce)$ at the point $q_0 = s/e = t/Ce$. Hence, if $q_0 \leq p$, that is, if $t \leq Cep$,

$$\mu(|f| \geq t) \leq 2e^{-\varphi(q_0)} = 2e^{-t/Ce}.$$

In case $q_0 \geq p$, that is, when $t \geq Cep$, just use (4.6) with (optimal) value $q = p$.

Similarly, under (4.2), if $2 \leq q \leq p$, we have $\mathbb{E}|f|^q \leq (C\sqrt{q})^q$. If $0 \leq q \leq 2$, use $\|f\|_q \leq C$ (by the weighted Poincaré) to get $\mathbb{E}|f|^q \leq 2(C\sqrt{q})^q$, where we have applied $2q^{q/2} \geq 2e^{-1/4e} > 1$. Hence, by Chebyshev's inequality,

$$\mu(|f| \geq t) \leq \frac{\mathbb{E}|f|^q}{t^q} \leq 2\left(\frac{C\sqrt{q}}{t}\right)^q = 2e^{-\varphi(q)/2}$$

with the same $\varphi$, corresponding to $s = t^2/s^2$. It remains to optimize over $q \in (0, p]$ and the proof is complete. $\square$

The interval $[t_0, t_1] \subset (0, +\infty)$, where the tail of $f$ admits an exponential bound, can be replaced by the whole half-axis under stronger moment hypotheses on the weight function. In particular, we have the following result:

COROLLARY 4.3. *Under (4.2), assume $\int e^{w^2/\alpha} \, d\mu \leq 2$ for some $\alpha > 0$. Given that $\|f\|_{\mathrm{Lip}} \leq 1$ and $\int f \, d\mu = 0$, for all $t > 0$,*

$$\mu(|f| \geq t) \leq 2 e^{-t/K\alpha}, \tag{4.7}$$

*where $K$ is a positive universal constant.*

Indeed, the moment assumption on $w$ is equivalent to $\|w\|_p \leq K_1 \sqrt{p}$ with an arbitrary $p \geq 1$, up to some constant $K_1$. Hence, by Theorem 4.1, $\|f\|_p \leq K_1 p$, which in turn is equivalent to (4.7) up to some constant $K$.

We now illustrate Theorem 4.1 and its Corollary 4.2 on the example of the generalized Cauchy distributions. To better see the role of the dimension, rescale the coordinates and consider the image $\tilde{\nu}_\beta$ of $\nu_\beta$ under the map $x \to \sqrt{2\beta - n} x$. The probability measure $\tilde{\nu}_\beta$ has density

$$\frac{d\tilde{\nu}_\beta(x)}{dx} = \frac{1}{Z} \left(1 + \frac{|x|^2}{2\beta - n}\right)^{-\beta},$$

up to some normalizing factor $Z$, so that when $\beta$ grows to infinity, these measures approach the standard Gaussian distribution $\gamma_n$ on $\mathbb{R}^n$. (Therefore, one may expect that some properties of $\tilde{\nu}_\beta$ with sufficiently large $\beta$ do not depend on the dimension, similarly to the case of the Gaussian measure.)

For the measure $\tilde{\nu}_\beta$, the weighted Poincaré-type and the weighted logarithmic Sobolev inequalities in Theorems 3.1 and 3.4 take the form

$$\mathrm{Var}_{\tilde{\nu}_\beta}(g) \leq C_\beta \frac{\beta - n/2}{\beta - 1} \int |\nabla g(x)|^2 \left(1 + \frac{|x|^2}{2\beta - n}\right) d\tilde{\nu}_\beta(x), \tag{4.8}$$

$$\mathrm{Ent}_{\tilde{\nu}_\beta}(g^2) \leq 2 \frac{\beta - n/2}{\beta - 1} \int |\nabla g(x)|^2 \left(1 + \frac{|x|^2}{2\beta - n}\right)^2 d\tilde{\nu}_\beta(x), \tag{4.9}$$

for every bounded smooth $g$ on $\mathbb{R}^n$, whenever $\beta \geq n$ ($\beta > 1$). Assume $\beta \geq n + 1$, so that $\beta \geq 2$. In that case, $C_\beta \leq (\sqrt{2} + \sqrt{3})^2 < 10$ and the constant in (4.8) is bounded by $10 \frac{\beta - n/2}{\beta - 1} < 15$. Hence, $\tilde{\nu}_\beta$ shares the weighted Poincaré-type inequality with the weight function $w(x)^2 = 15(1 + \frac{|x|^2}{2\beta - n})$. To bound its $L^p$-norms, we use a general elementary formula

$$\int_{\mathbb{R}^n} \left(1 + \frac{|x|^2}{r^2}\right)^{-a} dx = \frac{n \omega_n r^n}{2} \frac{\Gamma(n/2) \Gamma(a - n/2)}{\Gamma(a)}$$

with parameters $a > n/2$, $r > 0$ (where $\omega_n$ denotes the volume of the Euclidean unit ball in $\mathbb{R}^n$). Let $p = 2q$ with $q$ a positive integer. If $q < \alpha \equiv \beta - n/2$, then

$$\int_{\mathbb{R}^n} \left(1 + \frac{|x|^2}{2\beta - n}\right)^q d\tilde{\nu}_\beta(x) = \frac{(\beta - 1) \cdots (\beta - q)}{(\alpha - 1) \cdots (\alpha - q)} \leq \left(\frac{\beta - q}{\alpha - q}\right)^q. \tag{4.10}$$

If $q \leq \beta - \frac{3}{4}n$, then $\frac{\beta-q}{\alpha-q} \leq 3$, and (4.10) gives $\|w\|_p \leq \sqrt{45} < 7$. Hence, by the first part of Corollary 4.2 with the parameters $p = 2[\beta - \frac{3}{4}n]$ and $C = 7$, for any $f$ on $\mathbb{R}^n$ such that $\|f\|_{\mathrm{Lip}} \leq 1$ and $\int f \, d\tilde{\nu}_\beta = 0$,

$$\tilde{\nu}_\beta(|f| \geq t) \leq \begin{cases} 2e^{-t/7e}, & \text{if } 0 \leq t \leq t_1, \\ 2\left(\dfrac{7p}{t}\right)^p, & \text{if } t \geq t_1, \end{cases}$$

where $t_1 = 7e[\beta - \frac{3}{4}n]$. Note that, in view of the assumption $\beta \geq n + 1$, the values of $p$ and $t_1$ are as large as at least a factor of $n$.

Now consider the weight function $w(x)^2 = \frac{\beta - n/2}{\beta - 1}(1 + \frac{|x|^2}{2\beta - n})^2$ appearing in (4.9). Again, if $1 \leq q \leq \beta - \frac{3}{4}n$, from (4.10) with $q = p$ we get $\|w\|_p \leq 3\sqrt{3/2} < 4$. Hence, by involving the second part of Corollary 4.2 with the parameters $p = [\beta - \frac{3}{4}n]$ and $C = 4$ we arrive at the following conclusion:

COROLLARY 4.4. *If* $\beta \geq n + 1$, *for any function* $f$ *on* $\mathbb{R}^n$ *such that* $\|f\|_{\mathrm{Lip}} \leq 1$ *and* $\int f \, d\tilde{\nu}_\beta = 0$,

$$\tilde{\nu}_\beta(|f| \geq t) \leq \begin{cases} 2e^{-t^2/32e}, & \text{if } 0 \leq t < t_0, \\ 2e^{-t/7e}, & \text{if } t_0 \leq t \leq t_1, \\ 2\left(\dfrac{7p}{t}\right)^p, & \text{if } t \geq t_1, \end{cases}$$

*where* $p = 2[\beta - \frac{3}{4}n]$, $t_0 = 4\sqrt{e[\beta - \frac{3}{4}n]}$ *and* $t_1 = 7e[\beta - \frac{3}{4}n]$.

Thus, on an interval of length at least $\sqrt{n}$ (and even larger when $\beta$ increases), we obtain a Gaussian decay for the tail functions. At the expense of the numerical constants in the Gaussian and exponential bounds of Corollary 4.4, the assumption $\beta \geq n + 1$ may be weakened to $\beta \geq c(n + 1)$ with $\frac{1}{2} < c < 1$, and then the resulting intervals and bounds will involve also the parameter $c$, although they cannot be made uniform in $c$. If $\beta = \frac{n+1}{2} + o(n)$ and $\beta - \frac{n+1}{2} \to +\infty$, the previous arguments should work as well, but we are going to be led to more sophisticated estimates.

**5. Weighted Cheeger and Poincaré-type.** As the Cauchy distributions share the weighted Poincaré-type inequalities and the associated concentration phenomena, one may ask whether these properties extend to more general families of probability distributions. Of particular interest is the family of $\kappa$-concave measures in the hierarchy of convex measures as described in the Introduction. As a result, it turns out to be possible to choose, up to multiplicative constants, the common weight function $1 + |x|^2$ for these measures to satisfy a weighted Poincaré-type inequality. The conclusion will be reached by means of a stronger Cheeger-type isoperimetric inequality of independent interest.

Namely, let $\mu$ be a probability measure on $\mathbb{R}^n$ with density

$$(5.1) \qquad p(x) = V(x)^{-\beta}, \qquad \beta > n,$$

where $V$ is an arbitrary positive convex function on $\mathbb{R}^n$, possibly taking an infinite value. Denote by $r$ the unique positive number such that $\mu(|x| \leq r) = 2/3$, that is, the $\mu$-quantile of order $2/3$ for the Euclidean norm. Then we have:

THEOREM 5.1. *For some universal constant $C$, the measure $\mu$ satisfies the weighted Poincaré-type inequality*

$$(5.2) \qquad \mathrm{Var}_\mu(g) \leq C \frac{\beta - n + 1}{\beta - n} \int |\nabla g(x)|^2 (r^2 + |x|^2) \, d\mu(x)$$

*for all bounded smooth functions $g$ on $\mathbb{R}^n$.*

Note, however, that one may lose some information on the correct asymptotics in the constant in front of the integral. For example, when $\mu = \nu_\beta$ is the Cauchy distribution with parameter $\beta > n$, a factor of order $1/\beta$ is absent on the right-hand side of (5.2) with respect to (3.1).

As a main point in the proof of Theorem 5.1, recall, as discussed in the Introduction, that any measure $\mu$ with density of the form (5.1) satisfies the Brunn–Minkowski-type inequality

$$(5.3) \qquad \mu(tA + (1-t)B) \geq [t\mu(A)^\kappa + (1-t)\mu(B)^\kappa]^{1/\kappa}$$

for all $t \in (0,1)$ and for all Borel measurable sets $A, B \subset \mathbb{R}^n$ with $\kappa = -1/(\beta - n)$.

PROOF OF THEOREM 5.1. The first step consists in the stronger weighted Cheeger-type inequality

$$(5.4) \qquad \int |g| \, d\mu \leq D \int |\nabla g(x)|(r + |x|) \, d\mu(x)$$

with positive constants $D$ and $r$, where $g$ is an arbitrary smooth function on $\mathbb{R}^n$ with $\mu$-median zero. Splitting $g$ into the positive and negative parts, one may assume $g$ is nonnegative. Moreover, (5.4) is equivalent to its particular case where $g$ approximates the characteristic function of an arbitrary "nice" set $A \subset \mathbb{R}^n$ with measure $0 < \mu(A) \leq 1/2$ (cf. [24]). For example, it is enough to consider the sets $A$ that are representable as a finite union of closed balls, contained in the open supporting set of $\mu$. Then (5.4) reduces to the isoperimetric-type inequality

$$(5.5) \qquad \mu(A) \leq D\nu^+(A),$$

where $\nu^+(A)$ denotes the perimeter of the set $A$ with respect to the measure $d\nu(x) = (r + |x|) \, d\mu(x)$, namely

$$\nu^+(A) = \lim_{\varepsilon \downarrow 0} \frac{\nu(A + \varepsilon B_1) - \nu(A)}{\varepsilon} = \int_{\partial A} p(x)(r + |x|) \, d\mathcal{H}_{n-1}(x),$$

where $B_1$ is the unit Euclidean ball with centre at the origin and where $\mathcal{H}_{n-1}$ denotes the Lebesgue measure on the boundary $\partial A$ of the set $A$.

Introduce also the $\mu$-perimeter $\mu^+(A) = \int_{\partial A} p(x) \, d\mathcal{H}_{n-1}(x)$. Note that, for $\mathcal{H}_{n-1}$-almost all points $x$ in $\partial A$, the outer normal unit vector $n_A(x)$ at $x$ is well-defined, and for any $r > 0$,

$$\lim_{\varepsilon \downarrow 0} \frac{\mu((1-\varepsilon)A + \varepsilon B_r) - \mu(A)}{\varepsilon} = r\mu^+(A) - \int_{\partial A} \langle n_A(x), x \rangle p(x) \, d\mathcal{H}_{n-1}(x)$$

$$\leq \int_{\partial A} (r + |x|) p(x) \, d\mathcal{H}_{n-1}(x) = \nu^+(A).$$

On the other hand, the above left-hand side may be bounded by virtue of the geometric inequality (5.3). Applying the latter to the sets $A$ and $B = B_r$ (the Euclidean ball of radius $r$), we get that

$$\mu((1-\varepsilon)A + \varepsilon B_r) - \mu(A) \geq [(1-\varepsilon)\mu(A)^\kappa + \varepsilon\mu(B_r)^\kappa]^{1/\kappa} - \mu(A)$$

$$= -\frac{\varepsilon}{\kappa}\mu(A)^{1-\kappa}[\mu(A)^\kappa - \mu(B_r)^\kappa] + o(\varepsilon).$$

Letting $\varepsilon \to 0$ and combining the two inequalities, we get

(5.6) $$\nu^+(A) \geq \mu(A)^{1-\kappa} \frac{\mu(A)^\kappa - \mu(B_r)^\kappa}{-\kappa}.$$

Put $t = \mu(A)$, $p_r = \mu(B_r)$, and assume $0 < t \leq 1/2 < p_r \leq 1$. To adapt (5.6) to (5.5), we need a bound of the form $\frac{t^\kappa - p_r^\kappa}{-\kappa} \geq c_\kappa t^\kappa$ with some positive constant $c_\kappa$, which would give $\nu^+(A) \geq c_\kappa \mu(A)$. Clearly, the value $t = 1/2$ is critical, so for the optimal constant we have

$$c_k = \frac{1 - (2p_r)^\kappa}{-\kappa} \geq \frac{\log(2p_r)}{1 - \kappa}.$$

Hence (5.5) and (5.4) hold true with arbitrary $r > 0$, such that $\mu(B_r) > 1/2$, in which case we may put $D = (1 - \kappa)/(\log 2\mu(B_r))$.

We now conclude the proof and reach the weighted Poincaré-type inequality (5.2). To this task, if $g$ is an arbitrary nonnegative smooth function on $\mathbb{R}^n$ with $\mu$-median zero, apply (5.4) to $g^2$ to get with the help of the Cauchy-Schwarz inequality that

$$\int g^2 \, d\mu \leq 4D^2 \int |\nabla g(x)|^2 (r + |x|)^2 \, d\mu(x) \leq 8D^2 \int |\nabla g(x)|^2 (r^2 + |x|^2) \, d\mu(x).$$

At this step the assumption that $g$ is nonnegative may be removed, and when $\mu(B_r) \geq 2/3$, we arrive at the inequality (5.2) with constant $C = 8D^2 = 8/\log^2(\frac{4}{3})$. Theorem 5.1 is established in this way. □

The quantile $r$ of order $2/3$ in Theorem 5.1 may be replaced, up to numerical constants, with the median and other quantiles. It may also be replaced with integral characteristics, such as the moments $m_q = (\int |x|^q \, d\mu(x))^{1/q}$, $q > 0$, since by Chebyshev's inequality, $r \leq C m_q$ for some $C = C(q)$. For probability measures $\mu$ with densities (5.1), it was shown by Borell [8] that $\mu(|x| \geq t) = O(t^{1/\kappa})$ as $t \to +\infty$. This implies that $m_q$ are finite, whenever $q < \beta - n$. Moreover, as recently shown in [5], the constant $C$ can be made dependent on $\beta$, but independent of $q$. That is, we always have $r \leq C_\beta m_0$ where

$$m_0 = \lim_{q \to 0} m_q = \exp \int \log |x| \, d\mu(x)$$

is the geometric mean for the Euclidean norm under $\mu$.

COROLLARY 5.2. *Any convex probability measure $\mu$ on $\mathbb{R}^n$ with density (5.1) satisfies the weighted Poincaré-type inequality*

$$\mathrm{Var}_\mu(g) \leq C_\beta \int |\nabla g(x)|^2 (m_0^2 + |x|^2) \, d\mu(x)$$

*for all bounded smooth functions $g$ on $\mathbb{R}^n$, where the constant $C_\beta$ depends on $\beta$, only.*

In fact, with the help of a routine localization argument and a careful treatment of the one-dimensional case, this inequality may be sharpened for $\kappa = -1/(\beta - n)$ approaching zero as

(5.7) $$\mathrm{Var}_\mu(g) \leq C_\kappa \int |\nabla g(x)|^2 (m_0^2 + \kappa^2 |x|^2) \, d\mu(x)$$

with some $C_\kappa$ continuously depending on $\kappa \leq 0$. A slight advantage of this form is that the limit case $\kappa = 0$ involves the usual Poincaré-type inequality for an arbitrary log-concave probability measure $\mu$, namely

$$\mathrm{Var}_\mu(g) \leq C_0 m_0^2 \int |\nabla g(x)|^2 \, d\mu(x).$$

In an equivalent form, the latter was obtained in 1995 by Kannan, Lovász and Simonovits [15] by virtue of the localization lemma of [19], compare [4]. Another motivation to gain a small factor in front of $|x|^2$ in (5.7) is that we may then apply Proposition 3.3. The latter implies that, if $\kappa_0 < \kappa < 0$,

$$\inf_{c \in \mathbb{R}} \int \frac{|g(x) - c|^2}{m_0^2 + \kappa^2 |x|^2} \, d\mu(x) \leq C_\kappa \int |\nabla g|^2 \, d\mu,$$

where $\kappa_0$ is a numerical constant and $C_\kappa$ continuously depends on $\kappa$. It would be interesting to explore whether this property holds true for the whole range of $\kappa$.

SCHOOL OF MATHEMATICS
UNIVERSITY OF MINNESOTA
MINNEAPOLIS, MINNEAPOLIS 55455
USA
E-MAIL: bobkov@math.umn.edu

INSTITUT DE MATHÉMATIQUES
UNIVERSITÉ PAUL-SABATIER
31062 TOULOUSE
FRANCE
E-MAIL: ledoux@math.univ-toulouse.fr